\documentclass[11pt]{article}
\usepackage[english]{babel} %

\usepackage{tikz}
\usepackage{amssymb}
\usepackage{amsmath}
\usepackage{amsthm}
\usepackage{enumitem}
\usepackage{url}

\theoremstyle{plain}
\newtheorem{theorem}{Theorem}[section] % reset theorem numbering
\newtheorem*{theorem*}{Theorem} % reset theorem not numbering

\newtheorem{conjecture}[theorem]{Conjecture}
\newtheorem{lemma}[theorem]{Lemma}
\newtheorem{corollary}[theorem]{Corollary}

\theoremstyle{definition}
\newtheorem{remark}[theorem]{Remark}

\def\gr{\operatorname{gr}}

\def\GAP{\texttt{GAP}\ }

\title{Factorizations of simple groups of order 168 and 360}
\author{Mikhail Kabenyuk}
\date{}

\begin{document}\maketitle

\begin{abstract}
A finite group $G$ is called $k$-factorizable 
if for any factorization $|G|=a_1\cdots a_k$ with $a_i>1$
there exist subsets $A_i$ of $G$ with $|A_i|=a_i$ such that
$G=A_1\cdots A_k$.
We say that $G$ is \textit{multifold-factorizable}
if $G$ is $k$-factorizable for any possible integer $k\geq2$.

We prove that simple groups of orders 168 and 360 are
multifold-factorizable and formulate two conjectures
that the symmetric group $S_n$ for any $n$ and the alternative group $A_n$ 
for $n\geq6$ are multifold-factorizable.
\end{abstract}

\section{Introduction}\label{section-Introduction}
\def\Hajos{Haj$\acute{\rm{o}}$s}
\def\Redei{R$\acute{\rm{e}}$dei}
\def\Szabo{Szab$\acute{\rm{o}}$}

Let $G$ be a group
and let $A_1,\ldots,A_k$ be subsets of $G$ such that
each element $g\in G$ is uniquely represented as $g=a_1\ldots a_k$ with $a_i\in A_i$.
We write $G=A_1\ldots A_k$ and call this a \textit{$k$-factorization} of $G$.
If $G$ is a finite group and $n=|G|$, then 
to each $k$-factorization of the group $G$ there corresponds a $k$-factorization of $n$, i.e. 
if $G=A_1\ldots A_k$ and $m_i=|A_i|$, then $n=m_1\ldots m_k$.
In this case we will call this $k$-factorization $(m_1,\ldots,m_k)$-factorization.

Is the converse true?
We say that a finite group $G$ \textit{$k$-factorizable} 
if whatever factorization $|G|=m_1\ldots m_k$ with $m_i>1$ is, 
there exists a factorization of the group $G=A_1\ldots A_k$ with $m_i=|A_i|$.
If $G$ is $k$-factorizable for any possible integer $k\geq2$, then $G$ is called
\textit{multifold-factorizable}.

In answer to Hooshmand's question, 
(\cite[Problem 19.35]{Kourovka} see also \cite[Questions I and II]{Hooshmand}), 
Bergman \cite{Bergman} (see also \cite{Banakh}) proved that
the alternating group $A_4$ is not $3$-factorizable since it has no $(2,3,2)$-factorization.
By means of computer calculations, Brunault \cite{Banakh} proved that the alternating group $A_5$ does not have 
a $(2,3,5,2)$-factorization.
However, it does not yet follow that $A_5$ is not $3$-factorizable.

All finite groups are $2$-factorizable 
if every simple finite group is $2$-factorizable.
Each simple group of order at most $10000$ is
$2$-factorizable \cite{Bildanov}.

The author \cite{kabenyuk-1} proved the following statements: 
1) For each $k\geq3$ there exist of a finite group, which has no $k$-factorization. 
2) There are only $8$ non-multifold-factorizable groups 
among groups of order at most $100$.
3)
Let a finite group $G$ have the following properties:
$(i)$ a Sylow $2$-subgroup is elementary abelian;
$(ii)$ all involutions of $G$ are conjugate;
$(iii)$ the centralizer of every involution is the direct product of
a group of odd order and a Sylow $2$-subgroup of $G$.
Then $G$ is not $3$-factorizable.
It follows that all simple groups $SL(2,2^s)$, $s\geq4$,
are not $3$-factorizable (note that $SL(2,4)\cong A_5$).

In this context, it seems reasonable to ask the question: are there simple multifold-factorizable groups?

Here we prove that the simple groups of order $168$ and $360$ are
multifold-factorizable and formulate two conjectures
that the symmetric group $S_n$ for any $n\geq1$ and the alternative group $A_n$ 
for $n\geq6$ are multifold-factorizable.

The paper is organized as follows.
In the rest of this section, we will recall some of the notation.
In Section \ref{section-Preliminaries}, we will set out some of the necessary known facts.
In Sections \ref{section-group-168} and \ref{section-group-360} gives the proofs that
the simple groups of order $168$ and $360$ are multifold-factorizable.
In Section \ref{section-algorithms} we discuss how all factorizations of simple groups of order $168$ and $360$ were found.
Finally, in Section \ref{section-Questions} we formulate some conjectures on multifold-factorizable.

\phantom{a}

If $X$ is a subset of a group $G$, $|X|$ will denote the order of $X$
and $\gr(X)$ will denote the subgroup of $G$ generated by $X$.
If $X=\{x\}$, then $\gr(x)$ stands for cyclic subgroups of $G$ generated by $x$.
We write $H<G$ or $H\leq G$ whenever $H$ is a subgroup of $G$.
If $H<G$, then $|G:H|$ is the index of $H$ in $G$.
We denote an identity element in a group by the symbol $e$.

\section{Preliminaries}\label{section-Preliminaries}
We will need some results, which we will give here without detailed proofs.

\begin{lemma}\label{lemma_supersolvable_multifold}
    {\rm(\cite[Lemma 4.2]{kabenyuk-1}, \cite[Theorem 2.1]{Hooshmand})}
    If $G$ is a finite supersolvable group,
    then $G$ is multifold-factorizable.
\end{lemma}

\begin{lemma}\label{lemma_S4}
    {\rm(\cite[Lemma 5.4]{kabenyuk-1}, \cite[Example 2.2]{Hooshmand})}
    The symmetric groups $S_4$ is multifold-factorizable.
\end{lemma}

\begin{lemma}\label{lemma_ATB}{\rm(\cite[Lemma 5.6]{kabenyuk-1})}
    Let $G$ be a group of order $n$ and let $A$ and $B$ be subgroups of $G$
    of orders $a$ and $b$ respectively.
    Let
    $$
    G=At_1B\cup At_2B\cup\ldots\cup At_sB
    $$
    be a decomposition of $G$ into double cosets of $A$ and $B$,
    $T=\{t_1,t_2,\ldots,t_s\}$.
    Then $G=ATB$ is an $(a,s,b)$-factorization of $G$ if and only if
    $A^x\cap B=\{e\}$ for every $x\in G$.
\end{lemma}

\begin{corollary}\label{corollary_ATB}
    If $A$ and $B$ are subgroups of $G$
    of orders $a$ and $b$ respectively, $a$ and $b$ are relatively prime, and $s=\dfrac{|G|}{ab}$,
    then $G$ admits an $(a,s,b)$-factorization.
\end{corollary}

\begin{lemma}\label{lemma_Omega}
    Let $G$ be a group of order $n$ and
    let $m=\Omega(n)$ be
    the number of prime factors of $n$ (with multiplicity).
    If the group $G$ is $m$-factorizable, then $G$ is multifold-factorizable.
\end{lemma}

\begin{lemma}\label{lemma_SimpleGeneral}
Let $G$ be a finite group.
\begin{enumerate}[label=$(\roman*)$]
    \item\label{inverse_order}
If $G=A_1\cdot\ldots\cdot A_k$ is an $(m_1,\ldots,m_k)$-factorization of $G$, then
$G=A_k^{-1}\cdot\ldots\cdot A_1^{-1}$ is an $(m_k,\ldots,m_1)$-factorization.
    \item\label{subgroup_index}
If $H$ is a subgroup of $G$, $m=|G:H|$ and $H$ has
an $(m_1,\ldots,m_k)$-factorization, then
the group $G$ possesses factorizations of the forms
$(m_1,\ldots,m_k,m)$ and $(m,m_1,\ldots,m_k)$.
    \item\label{normal_subgroup}
If $H$ is a normal subgroup of $G$ and the factor group $G/H$
has an $(m_1,\ldots,m_k)$-factorization, then
the group $G$ has factorizations of each of the following forms
$$
(|H|,m_1,\ldots,m_k),\
(m_1,|H|,m_2,\ldots,m_k),\
\ldots,
(m_1,\ldots,m_k,|H|).
$$
\end{enumerate}
\end{lemma}

\begin{corollary}\label{corollary_prime_index}
Let $G$ be a finite group and $H$ be a subgroup of prime index $p$.
If $H$ is multifold-factorizable, then $G$ has any factorization
except maybe factorizations of the form $(p_1,*,\ldots,*,p_2)$,
where primes $p_1$ and $p_2$ are distinct from $p$.
\end{corollary}

\begin{corollary}\label{corollary_normal_prime_order}
Let $G$ be a finite group and $H$ be a normal subgroup of prime order $p$.
If the factor group $G/H$ is multifold-factorizable, then $G$ is multifold-factorizable too.
\end{corollary}

The following statements are useful for the reduction of computations in the search for factorizations.
\begin{lemma}\label{lemma_divisibility}
    {\rm(\cite[Lemma 1.2]{Bergman}, \cite[Lemma 2.1]{kabenyuk-1})}    
    If $G=A_1\ldots A_k$ is a factorization of a finite group $G$, then
    the orders of subgroups $\gr(A_1)$ and $\gr(A_k)$ are divisible by $|A_1|$ and $|A_k|$, respectively.
\end{lemma}
\textit{Proof}.
If $X=\{x\in A_2\ldots A_k\mid A_1x\subset\gr(A_1)\}$, then
$$
\cup_{x\in X}A_1x=\gr(A_1).
$$
If  $x,y\in A_2\ldots A_k$ and $x\neq y$, then $A_1x\cap A_1y=\varnothing$.
Hence $|\gr(A_1)|$ is a multiple of $|A_1|$.
The statement about $\gr(A_k)$ is seen in the same way.

\begin{corollary}\label{corollary_on two element factors}
    If $G=A_1\ldots A_k$ is a factorization of a finite group $G$ and $A_1=\{e,a\}$
    where $e$ is the identity element of $G$, then
    $\gr(a)$ is of even order. The same statement is true for $A_k$.
\end{corollary}

\begin{lemma}\label{lemma_normalization}
    {\rm(\cite[Lemma 1.3]{Bergman}, \cite[Lemma 3.1]{kabenyuk-1})}    
    If a finite group $G$ has an $(m_1,\ldots,m_k)$-factorization, 
    then it has such an $(m_1,\ldots,m_k)$-factorization in which 
    each factor contains the identity element $e$.
\end{lemma}

Suppose we have $s$ of distinct letters $a_1,\ldots,a_s$
and $s$ of positive integers $k_1,\ldots,k_s$.  Let $n=k_1+k_2+\ldots+k_s$.
Let $W$ be the set of all those words in the alphabet $a_1,\ldots,a_s$ in which the letter $a_i$ occurs exactly $k_i$ times for each $1\leq i\leq s$.
It is well known that
\begin{equation}\label{formula_W}
  |W|=\frac{n!}{k_1!\ldots k_s!}.
\end{equation}
A word $w\in W$ is called a palindrome if it reads the same backward as forward.
Denote by $P$ the set of all palindromes in $W$.
We have

\begin{equation}\label{palindrome}
|P|=
\left\{
  \begin{array}{ll}
   \dfrac{(n/2)!}{(k_1/2)!\ldots(k_s/2)!},
        & \hbox{if all $k_i$ are even;} \\
   \dfrac{((n-1)/2)!}{(k_1/2)!\ldots((k_r-1)/2)!\ldots(k_s/2)!},
        & \hbox{if only $k_r$ is odd;} \\
    0, & \hbox{else}.
  \end{array}
\right.
\end{equation}
We will denote by $w'$ the word $w$ read backwards (if $w=x_1\ldots x_n$, then $w'=x_n\ldots x_1$).
Two words $w_1,w_2\in W$ are called equivalent if $w_2=w_1'$.
A word $w\in W$ is a palindrome if and only if $w'=w$.
Under these notations, we have
\begin{lemma}\label{lemma_number_equivalence_classes}
    The number of equivalence classes of words from $W$ is equal
    \begin{equation}
        \frac{|W|+|P|}{2},
    \end{equation}
    where $|W|$ and $|P|$ are given by formulas $(\ref{formula_W})$ and $(\ref{palindrome})$, respectively.
\end{lemma}

\section{Simple group of order $168$}\label{section-group-168}
As is well-known,
there is exactly one simple group of order $168$.
We will denote it here by $G$. There are different representations of this group,
for example $G\cong PSL(2,7)\cong GL(3,2)$.
The group $G$ is characterized
as the second smallest non-commutative simple group.
Here we represent $G$ as a permutation group on $7$ letters:
$$
G=\gr((3,4)(5,6),\ (1,2,3)(4,5,7)).
$$
This is how this group is defined in \texttt{GAP} with the command
\medskip
\newline
\hspace*{\parindent}\texttt{G:=AllSmallGroups(168,IsSimple,true)[1];}
\medskip
\newline
Note that the group $G$ contains a subgroup of order $24$ isomorphic to $S_4$,
and a non-abelian subgroup of order $21$ isomorphic to the semidirect product $C_7\rtimes C_3$.

\begin{theorem}
  The simple group $G$ of order $168$ is multifold-factorizable.
\end{theorem}
\textit{Proof.}
Since $\Omega(168)=5$,
it suffices to show that the group $G$ is $5$-factorizable (Lemma \ref{lemma_Omega}).
By Lemma \ref{lemma_number_equivalence_classes} the number of distinct representations of $168$
as a product of primes is equal to $10$.
Since $G$ contains a subgroup $H$ isomorphic to $S_4$ and $S_4$ is
multifold-factorizable (Lemma \ref{lemma_S4}) and $H$ has index $7$ in $G$, it follows by Corollary \ref{corollary_prime_index} that
we can discard factorizations of $168$ that start or end with $7$.
Therefore, we only need to consider the following $6$ factorizations of $168$:
\begin{align*}
(2,2,2,7,3),&\quad
(2,2,3,7,2),\\
(2,3,2,7,2),&\quad
(3,2,2,7,2),\\
(2,3,7,2,2),&\quad
(3,2,7,2,2).
\end{align*}
The theorem will be completely proved if we construct for each of these factorisations the corresponding factorisation of the group $G$.
This will be done in the Lemmas \ref{lemma_22273}, \ref{lemma_6-7-2-2}, and \ref{lemma_6-2-7-2} below.

\begin{lemma}\label{lemma_22273}
    The group $G$ has a $(2,2,2,7,3)$-factorization.
\end{lemma}
\textit{Proof.}
Let $P$ be a Sylow $2$-subgroup of $G$, and let $Q$ be a subgroup of $G$ of order $21$.
We have $G=PQ$.
Both groups $P$ and $Q$ are supersolvable, by Lemma \ref{lemma_supersolvable_multifold} they are multifold-factorizable, and then the group $G$ has a $(2,2,2,7,3)$-factorization.
The lemma is proved.

For the remaining five cases, we do not know simple ways to compute the corresponding factorizations of the group $G$.
\begin{lemma}\label{lemma_6-7-2-2}
    The group $G$ has $(2,3,7,2,2)$- and $(3,2,7,2,2)$-factorizations.
\end{lemma}
\textit{Proof.}
Let $A=\gr(u,v)$, where $u=(2,3,4)(5,6,7)$ and $v=(3,4)(5,6)$.
It is easy to see that the group $A$ has order $6$.
Let $d=(2,7)(3,4,5,6)$ and $c=(1,6,3,5,2,4,7)$
and let $D=\{e,d\}$ and $C=\{e,c\}$.
Let finally
\begin{equation*}
    \begin{split}
    B&=\{e, (3,5)(4,6), (1,2)(3,6), (1,2,3)(4,5,7), \\
    &\hspace*{22pt} (1,2,5)(3,7,6), (1,5,3)(2,4,6), (1,6,3)(4,7,5)\}.
    \end{split}
\end{equation*}
Now we can check that $G=ABCD$.
It is clear that checking this equality by hand would be too tedious.
So we'll point out how we can do this check using \texttt{GAP} \cite{GAP}:
\medskip
\newline
\texttt{\small
\hspace*{\parindent}G:=AllSmallGroups(168,IsSimple,true)[1];\\
\hspace*{\parindent}A:=Group([(2,3,4)(5,6,7),(3,4)(5,6)]);;\\
\hspace*{\parindent}B:=[(), (3,5)(4,6), (1,2)(3,6), (1,2,3)(4,5,7), \\
\hspace*{\parindent}\qquad\ (1,2,5)(3,7,6), (1,5,3)(2,4,6), (1,6,3)(4,7,5)];;\\
\hspace*{\parindent}C:=[(), (1,6,3,5,2,4,7)];;\\
\hspace*{\parindent}D:=[(), (2,7)(3,4,5,6)];;\\
\hspace*{\parindent}P:=List(Cartesian(A,B,C,D),x-> Product(x));;\\
\hspace*{\parindent}Size(AsSet(P))=168;\\
\hspace*{\parindent}IsSubset(P,G);
}
\medskip
\newline
As a result of the above commands, the \texttt{GAP} will output $|P|=|ABCD|=168$, which means that the equality
$G=ABCD$  is true.
Since $|B|=7$, $|C|=|D|=2$, $A$ has $(2,3)$- and $(3,2)$-factorizations,
it follows that
$G$ has $(2,3,7,2,2)$- and $(3,2,7,2,2)$-factorizations and
the lemma is proved.

\begin{remark}\label{remark-OnGAP}
Generally speaking, we should make sure that all the elements
that appear in the proof of this lemma and the other lemmas of this section below
lie in the group $G$, i.e.
$$
A\leq G,\ c,d\in G,\ B\subset G.
$$
To check this we can use a \GAP command:
\texttt{IsSubset(P,G);}
The output should be \texttt{true}. 

Note also that our subgroup $A$ is the first member of the list of subgroups of order $6$ in the group $G$, that is, 
$$
\texttt{A:=First(AllSubgroups(G),x->Order(x)=6);}
$$
\end{remark}

\begin{lemma}\label{lemma_12-7-2}
    The group $G$ has $(2,2,3,7,2)$- and $(3,2,2,7,2)$-factorizations.
\end{lemma}
\textit{Proof.}
Let's choose a subgroup $A=\gr((1,2,7)(3,6,4),(3,4)(5,6))$ of order $12$ in the group $G$.
In addition, put $C=\{(),(2,4,7,6)(3,5)\}$ and
\begin{equation*}
    \begin{split}
   B&=\{(), (2,3)(5,7), (2,7)(4,6), (1,2,3)(4,5,7), \\
    &\hspace*{22pt} (1,2,4,6,5,7,3), (1,3)(2,6), (1,3)(2,7,6,4) \}.
    \end{split}
\end{equation*}
Let us use \GAP to compute the product $ABC$:
\medskip
\newline
\texttt{\small
\hspace*{\parindent}G:=AllSmallGroups(168,IsSimple,true)[1];\\
\hspace*{\parindent}A:=Group([(1,2,7)(3,6,4), (3,4)(5,6)]);;\\
\hspace*{\parindent}B:=[(), (2,3)(5,7), (2,7)(4,6), (1,2,3)(4,5,7), \\
\hspace*{\parindent}\qquad\ (1,2,4,6,5,7,3), (1,3)(2,6), (1,3)(2,7,6,4)];;\\
\hspace*{\parindent}C:=[(), (2,4,7,6)(3,5)];;\\
\hspace*{\parindent}P:=List(Cartesian(A,B,C),x-> Product(x));;\\
\hspace*{\parindent}Size(AsSet(P))=168;\\
\hspace*{\parindent}IsSubset(P,G);
}
\medskip
\newline
As a result, we get $|ABC|=168$, hence $G=ABC$.
Since the group $A$ has $(2,2,3)$- and $(3,2,2)$-factorizations, and $|B|=7$, it follows that
$G$ possesses $(2,2,3,7,2)$- and $(3,2,2,7,2)$-factorizations and
this proves the lemma.

\begin{remark}\label{remark-to lemma-14-7-2}
The subgroup $A$ of Lemma \ref{lemma_12-7-2} is the first member of the list of subgroups of order $12$ in the group $G$; that is, 
$$
\texttt{A:=First(AllSubgroups(G),x->Order(x)=12);}
$$ 
\end{remark}

\begin{lemma}\label{lemma_6-2-7-2}
    The group $G$ has $(2,3,2,7,2)$-factorization.
\end{lemma}
\textit{Proof.}
Let $A$ be the same group as in Lemma \ref{lemma_6-7-2-2},
$b=(1,5)(6,7)$, $d=(1,5)(2,7,4,6)$, and
\begin{equation*}
    \begin{split}
    C&=\{(), (3,5)(4,6), (1,4,5,2)(6,7), (1,4,6)(3,5,7), \\
    &\hspace*{22pt} (1,4,6)(2,7,5), (1,2,5)(3,7,6), (1,6,2,4,5,3,7)\}.
    \end{split}
\end{equation*}
Let $B=\{e,b\}$, $D=\{e,d\}$.
We specify the \GAP commands for checking the equality $G=ABCD$.
\medskip
\newline
\texttt{\small
\hspace*{\parindent}G:=AllSmallGroups(168,IsSimple,true)[1];\\
\hspace*{\parindent}A:=Group([(2,3,4)(5,6,7), (3,4)(5,6)]);;\\
\hspace*{\parindent}B:=[(), (1,5)(6,7)];;\\
\hspace*{\parindent}C:=[(), (3,5)(4,6), (1,4,5,2)(6,7), (1,4,6)(3,5,7),\\
\hspace*{\parindent}\qquad\ (1,4,6)(2,7,5), (1,2,5)(3,7,6), (1,6,2,4,5,3,7)];;\\
\hspace*{\parindent}D:=[(), (1,5)(2,7,4,6)];;\\
\hspace*{\parindent}P:=List(Cartesian(A,B,C,D),x-> Product(x));;\\
\hspace*{\parindent}Size(AsSet(P))=168;\\
\hspace*{\parindent}IsSubset(P,G);
}
\medskip 
\newline
It follows that $G$ has $(2,3,2,7,2)$-factorization and 
the lemma is proved.

\section{Simple group of order $360$}\label{section-group-360}
We know that $A_6$ is the only simple group of order $360$.
The group $A_6$ contains a subgroup of order $24$, isomorphic to $S_4$;
a non-abelian subgroup of order $36$ with trivial center, isomorphic to the semidirect product
$(C_3\times C_3)\rtimes C_4$;
and a subgroup of order $60$, isomorphic to $A_5$.

\begin{theorem}\label{theorem-360}
  The group $A_6$ is multifold-factorizable.
\end{theorem}
\textit{Proof.}
Since $\Omega(360)=6$,
it suffices to show that $A_6$ is $6$-factorizable (Lemma \ref{lemma_Omega}).
By Lemma \ref{lemma_number_equivalence_classes} the number of distinct representations of $360$
as a product of primes is equal $30$. 
All such representations can be divided into three parts of 10 each:
\begin{equation}\label{three parts}
  (*,*,*,*,*,5),\ (*,*,*,*,5,*),\ (*,*,*,5,*,*).   
\end{equation}
Here and below under the asterisks are hidden numbers $2$ and $3$ in the required quantity.

\begin{table}[htp]
\begin{center}
\begin{tabular}{|r|c|c|c|c|}
  \hline
       & Factorization of $360$ & $A$   & $B$ & Lemma\\
  \hline
   1-4  &  $(*,*,*,*,{\mathbf 3},5)$  &$S_4$  &$C_5$  &\\
   5-10 &  $(*,*,*,*,{\mathbf 2},5)$  &$H$  &$C_5$  &\\
   11-14&  $(*,*,*,*,{\mathbf 5},3)$  &$S_4$  &$C_3$  &\\
   15   &  $(2,2,3,3,5,2)$  &$ $  & $ $  & Lemma \ref{lemma_A223352-322352}\\
   16   &  $(2,3,2,3,5,2)$  &$ $  & $ $  & Lemma \ref{lemma_A232352}\\
   17   &  $(3,2,2,3,5,2)$  &$ $  & $ $  & Lemma \ref{lemma_A223352-322352}\\
   18   &  $(2,{\mathbf 3},3,2,5,2)$  &$\{e,v\}$  & $A_5,\ (3,\ldots)$  &Lemma \ref{lemma_A5-[e,v]}\\
   19   &  $(3,{\mathbf 2},3,2,5,2)$  &$C_3 $  & $A_5,\ (3,\ldots)$  &\\       
   20   &  $(3,3,2,2,5,2)$  &$ $  & $ $  &Lemma \ref{lemma_A332252}\\
   21   &  $(2,2,2,{\mathbf 5},3,3)$  &$P_2$  &$P_3$  &\\   
   22-24&  $(*,*,*,5,{\mathbf 2},3)$  &$A_5,\ (\ldots,5)$  &$C_3$  &\\
   25-27&  $(*,*,*,5,{\mathbf 3},2)$  &$A_5,\ (\ldots,5)$  &$\{e,v\}$  &Lemma \ref{lemma_A5-[e,v]}\\
   28-30&  $(*,*,*,{\mathbf 5},{\mathbf 2},2)$  &$F$  &$\{e,v\}$  &Lemma \ref{lemma_A18-522}\\   
  \hline
\end{tabular}
\caption{6-factorizable $A_6$}
\label{Factorizable}
\end{center}
\end{table}
We use Lemma \ref{lemma_ATB} or Corollary \ref{corollary_ATB}.
It is sufficient to choose appropriate subgroups $A$ and $B$ of $A_6$.
This choice is specified in Table \ref{Factorizable}, where
we take the following notations:
$P_q$ is a Sylow $q$-subgroup of $A_6$, $q=2,3$;
$S_4=\gr((1,2,3),(1,2,3,4)(5,6))$ is isomorphic to the symmetric group of degree $4$;
$A_5$ is the alternating group on $\{1,\ldots,5\}$;
$H$ is a subgroup of order $36$, e.g. $H=\gr((1,2,3),(4,5,6),(1,4,2,5)(3,6))$;
$F$ is a subgroup of order $18$, e.g. $F=\gr((1,2,3),(4,5,6),(1,2)(4,5))$;
$C_3=\gr((1,2,3)(4,5,6))$;
$C_5=\gr(1,2,3,4,5)$.
In all cases it is clear that $A^x\cap B=\{e\}$ for every $x\in A_6$.
Number in boldface type in the second column of this table
is the number of double cosets of $AtB$ in $A_6$.

We know that all $p$-groups and the groups of order $24$ are multifold-factorizable 
\cite[Lemma 5.4]{kabenyuk-1}.
Due to \cite[Theorem 1.3]{kabenyuk-1}, we know that every group of order $36$ except two is factorizable. 
Both exceptional groups of order $36$ have a nontrivial center. 
The center of our group $H$ is trivial. 
Hence by \cite[Theorem 1.3]{kabenyuk-1} $H$ is multifold-factorizable.
All groups of order $18$ are also multifold-factorizable since they are all supersolvable, 
in particular our group $F$ is multifold-factorizable. 

At the same time, the group $A_5$ is not multifold-factorizable \cite[Theorem 1.3]{kabenyuk-1}. However, there are only three factorizations of $(2,15,2)$, $(2,3,5,2)$, and $(2,5,3,2)$ of number $60$ for which there is no corresponding factorization of the group $A_5$ (see Lemma \ref{lemma_A5}).
Therefore, the table \ref{Factorizable}, where the group $A_5$ acts as $A$ or $B$, also briefly indicates which factorization of the group $A_5$ is used.

In three cases a subset of two elements $\{e,v\}$ plays the role of $A$ or $B$, 
here $e$ is the identity permutation and $v=(1,2,3,6)(4,5)$.
In these cases, our problem is solved by Lemmas \ref{lemma_A5-[e,v]} and \ref{lemma_A18-522} formulated below.

However, there are 4 more cases in Table 1, these are rows 15, 16, 17, and 20 where groups $A$ and $B$ are not specified. These are the cases which cannot be solved by the lemma on double cosets. 
These 4 cases are discussed in lemmas \ref{lemma_A332252} -- \ref{lemma_A232352}.
Thus by proving these lemmas we obtain a complete proof of Theorem \ref{theorem-360}.

\begin{lemma}\label{lemma_A5}
    For any factorization $(m_1,\ldots,m_k)$ of number $60$ different 
    from $(2,15,2)$, $(2,3,5,2)$, and $(2,5,3,2)$ there exists 
     a factorization of the group $A_5=X_1\ldots X_k$ with $|X_i|=m_i$. 
\end{lemma}
\textit{Proof}. 
By virtue of lemma \ref{lemma_number_equivalence_classes}, 
it suffices to construct factorizations of the group $A_5$ for the following factorizations of the number $60$:
$$
(2,2,3,5),\
(2,3,2,5),\
(3,2,2,5),\
(2,2,5,3),\
(3,2,5,2),\
$$
Since $A_5=A_4\cdot C_5$ and $A_4$ has $(2,2,3)$- and $(3,2,2)$-factorizations,
it follows that $A_5$ admits $(2,2,3,5)$- and $(3,2,2,5)$-factorizations.
Let $P$ be a Sylow $2$-subgroup and $H$ be a subgroup of order $10$ of $A_5$.
By Corollary \ref{corollary_ATB} the number $(P,C_3)$-double cosets is $5$ and  
the number $(H,C_3)$-double cosets is $2$,
therefore there exist $(2,2,5,3)$- and $(3,2,5,2)$-factorizations of $A_5$.

Let us proceed to construct a $(2,3,2,5)$-factorization.
It seems natural to take this approach. 
Choose a subgroup $H$ of order $10$ in $A_5$ and then try to construct 
two sets $A$ and $B$ of size $2$ and $3$ respectively so that $A_5=ABH$. 
However, this way is doomed to failure, because in this case 
we will get also $(2,3,5,2)$-factorization of $A_5$, which as we know is impossible.

Therefore, we are forced to take another way. 
Let $a=(1,2)(3,4)$, $A=\gr(a)$, and $d=(1,2,3,4,5)$, $D=\gr(d)$.
By Corollary \ref{corollary_ATB} the number $(A,D)$-double cosets is $6$
and here's a list of representatives of all the $(A,D)$-double cosets:
$$
e, (3,4,5), (3,5,4), (2,3)(4,5), (2,3,4), (2,3,5). 
$$
Thus, if $b=(3,4,5)$, $c=(2,3)(4,5)$, and $B=\gr(b)$, $C=\gr(c)$, 
then $A_5=ABCD$ is a $(2,3,2,5)$-factorization of group $A_5$.
The last equality can be checked using the following \GAP commands:
\medskip
\newline
\texttt{\small
\hspace*{\parindent}A:=Group([(1,2)(3,4)]);;\\
\hspace*{\parindent}B:=Group([(3,4,5)]);;\\
\hspace*{\parindent}C:=Group([(2,3)(4,5)]);;\\
\hspace*{\parindent}D:=Group([(1,2,3,4,5)]);;\\
\hspace*{\parindent}P:=List(Cartesian(A,B,C,D),x-> Product(x));;\\
\hspace*{\parindent}Size(AsSet(P))=60;
}
\medskip
\newline 
We see that in this case all factorization multipliers are subgroups and
the lemma is proved.

\begin{lemma}\label{lemma_A5-[e,v]}
  If $v=(1,2,3,6)(4,5)$, $t=(1,2,3)(4,5,6)$ and $A_5$ is the alternating group on $\{1,\ldots,5\}$, then 
  
\begin{equation}\label{eq-60-3-2}
  A_6=A_5\cdot\{e,v^2,t\}\cdot\{e,v\}=\{e,v\}\cdot\{e,v^2,t\}\cdot A_5.
\end{equation}
\end{lemma}
\textit{Proof}. 
We have 
\begin{align*}
\{e,v^2,t\}\cdot\{e,v\}=\{e,v,v^2,v^3,t,tv\}=T_1,\\
\{e,v\}\cdot\{e,v^2,t\}=\{e,v,v^2,v^3,t,vt\}=T_2.
\end{align*}
It is easy to see that for every $i$, $1\leq i\leq6$, there exists $\sigma\in T_1$ such that $\sigma(6)=i$ and there exists $\tau\in T_2$ such that $\tau(i)=6$.
It follows that $T_1$ is a right transversal and
$T_2$ is a left transversal to $A_5$ in $A_6$ and 
Lemma is proved.
\begin{remark}
We can also check for equality (\ref{eq-60-3-2}) by the following \GAP commands:
\medskip
\newline
\texttt{\small
\hspace*{\parindent}A:=AlternatingGroup(5);;\\
\hspace*{\parindent}v:=(1,2,3,6)(4,5);;\\
\hspace*{\parindent}t:=(1,2,3)(4,5,6);;\\
\hspace*{\parindent}B:=[ (), v\textasciicircum2, t ];;\\
\hspace*{\parindent}C:=[ (), v ];;\\
\hspace*{\parindent}P:=List(Cartesian(A,B,C),x-> Product(x));;\\
\hspace*{\parindent}Q:=List(Cartesian(C,B,A),x-> Product(x));;\\
\hspace*{\parindent}Size(AsSet(P))=360;Size(AsSet(Q))=360;
}
\medskip
\end{remark}

\begin{lemma}\label{lemma_A18-522}
The group $A_6$ admits an $(18,5,2,2)$-factorization with the first factor equal to 
the subgroup $F=\gr((1,2,3), (4,5,6), (1,2)(4,5))$ of order $18$.
\end{lemma}
\textit{Proof.}
Let $v=(1,2,3,6)(4,5)$, 
$$
B=\{e, (3,4)(5,6), (3,5,4), (2,3,4), (2,4)(3,5)\},
$$
 and
$C=\{e, (1,3)(2,5,4,6)\}$.
Let us also denote by $D=\{e,v\}$.
Using the \GAP commands we check that $A_6=FBCD$:
\medskip
\newline
\texttt{\small
\hspace*{\parindent}F:=Group([(1,2,3), (4,5,6), (1,2)(4,5)]);;\\
\hspace*{\parindent}B:=[ (), (3,4)(5,6), (3,5,4), (2,3,4), (2,4)(3,5) ];;\\
\hspace*{\parindent}C:=[ (), (1,3)(2,5,4,6) ];;\\
\hspace*{\parindent}D:=[ (), (1,2,3,6)(4,5) ];;\\
\hspace*{\parindent}P:=List(Cartesian(F,B,C,D),x-> Product(x));;\\
\hspace*{\parindent}Size(AsSet(P))=360;
}
\medskip
\newline
Hence, $A_6$ possesses an $(18,5,2,2)$-factorization and 
the lemma is proved.

\begin{lemma}\label{lemma_A332252}
    The group $A_6$ admits a $(3,3,2,2,5,2)$-factorization.
\end{lemma}
\textit{Proof.}
Let $A=\gr((1,2,3),(4,5,6))$, $D=\gr((1,2,3,4,5),(1,5)(2,4))$, 
$B=\{e,(2,3,4\}$, $C=\{e,(3,6,4)\}$, .
We have $|A|=9$, $|D|=10$, and $|B|=|C|=2$.
Check that $A_6=ABCD$:
\medskip
\newline
\texttt{\small
\hspace*{\parindent}A:=Group([(1,2,3),(4,5,6)]);;\\
\hspace*{\parindent}B:=[ (), (2,3,4) ];;\\
\hspace*{\parindent}C:=[ (), (3,6,4) ];;\\
\hspace*{\parindent}D:=Group([(1,2,3,4,5),(1,5)(2,4)]);;\\
\hspace*{\parindent}P:=List(Cartesian(A,B,C,D),x-> Product(x));;\\
\hspace*{\parindent}Size(AsSet(P))=360;
}
\medskip
\newline
It follows that $A_6$ has a $(3,3,2,2,5,2)$-factorization and 
the lemma is proved.

\begin{lemma}\label{lemma_A223352-322352}
    The group $A_6$ admits a $(2,2,3,3,5,2)$-factorization and
a $(3,2,2,3,5,2)$-factorization.
\end{lemma}
\textit{Proof}.
Let 
\begin{align*}
&A=\gr((1,2,3),(1,2,4)),\\ 
&B=\{e, (3,4,6), (1,2)(3,4,6,5)\},\\
&C=\{e, (2,4,6,3,5), (1,6,5,4,3), (1,6,2,5,3), (1,6,4,2,3)\},\\
&D=\{e,(1,2,3,4)(5,6)\}.
\end{align*}
The group $A$ has order $12$, $|B|=3$, $|C|=5$, and $|D|=2$.
The following \GAP commands allow you to check the equality $A_6=ABCD$:
\medskip
\newline
\texttt{\small
\hspace*{\parindent}A:=Group([(1,2,3),(1,2,4)]);;\\
\hspace*{\parindent}B:=[ (), (3,4,6), (1,2)(3,4,6,5) ];;\\
\hspace*{\parindent}C:=[ (), (2,4,6,3,5), (1,6,5,4,3), (1,6,2,5,3), (1,6,4,2,3) ];;\\
\hspace*{\parindent}D:=[ (), (1,2,3,4)(5,6)];;\\
\hspace*{\parindent}P:=List(Cartesian(A,B,C,D),x-> Product(x));;\\
\hspace*{\parindent}Size(AsSet(P))=360;
}
\medskip
\newline
After executing all the above commands, the output is $360$.
This means that $|ABCD|=360$ and therefore $A_6=ABCD$.
Since the group $A$ admits $(2,2,3)$- and $(3,2,2)$-factorizations, 
it follows that $A_6$ has the specified factorizations and
the lemma is proved.

\begin{lemma}\label{lemma_A232352}
    The group $A_6$ admits a $(2,3,2,3,5,2)$-factorization.
\end{lemma}
\textit{Proof}.
Let $A=\{e,t,t^2,a,at,at^2\}$, where $a=(1,2)(3,4,5,6)$, $t=(1,3,5)$, and
$B=\{e,(1,2,6)\}$, 
$C=\{e,(3,4,6),(2,3,4)\}$
and
$$
D=\gr((1,2)(3,4),(1,3,5,4,2)).
$$
To check the equality $A_6=ABCD$, we use the \GAP commands:
\medskip
\newline
\texttt{\small
\hspace*{\parindent}a:=(1,2)(3,4,5,6);;\\
\hspace*{\parindent}t:=(1,3,5);;\\
\hspace*{\parindent}A:=[ (),t,t\textasciicircum2,a,a*t,a*t\textasciicircum2];;\\
\hspace*{\parindent}B:=[ (), (1,2,6) ];;\\
\hspace*{\parindent}C:=[ (), (3,4,6), (2,4,3) ];;\\
\hspace*{\parindent}D:=Group([(1,2)(3,4),(1,3,5,4,2)]);;\\
\hspace*{\parindent}P:=List(Cartesian(A,B,C,D),x-> Product(x));;\\
\hspace*{\parindent}Size(AsSet(P))=360;
}
\medskip
\newline
Therefore, $A_6$ possesses a $(6,2,3,10)$-factorization.
Since $A=\{e,a\}\cdot\{e,t,t^2\}$ 
and $D$ is a multifold-factorizable group of order $10$, 
it follows that $A_6$ admits a $(2,3,2,3,5,2)$-factorization
and the lemma is proved.

\section{Methods for the factorization of groups}\label{section-algorithms}
I note at once that I do not have a universal algorithm 
for factorization of an arbitrary finite group. 
The purpose of this section is to describe how all factorizations 
of simple groups of order $168$ and $360$ were found. 
Here, in addition to methods relying on the knowledge of subgroups, a brute force search was used. 
In those cases when the brute force search led to computations which my modest computer could not cope with in reasonable time, I used a trick from graph theory.

First of all, we always look for a factorization in which all the factors contain an identity element.
Next, we proceed from the fact that if a group $G$ has a subgroup $H$ and this subgroup has an 
$(m_1,\ldots,m_s)$-factorization, then it is very plausible that the group $G$ has an $(m_1,\ldots,m_s,m_{s+1},\ldots,m_k)$-factorization and also 
an $(m_{s+1},\ldots,m_k,m_1,\ldots,m_s)$-factorization for some $m_{s+1},\ldots,m_k$.

Due to this observation, we manage in many cases to reduce the whole computation 
to finding factorizations $G=ABCD$ consisting of four factors, 
and as a rule either $A$ or $D$, or $A$ and $D$ are both subgroups of $G$. 
In any case, one always manages to guess what $A$ and $D$ are even if one of them is not a subgroup.

\textsc{Case 1.}
Let $A$ and $D$ be subgroups of $G$ with coprime orders, 
or at least the condition 
$
xAx^{-1}\cap D=\{e\}
$
for each $x\in G$ is satisfied.
Let $W$ be a right transversal of $A$ in $G$ and $V$ be a left transversal of $D$ in $G$.
Let further $W'$ (resp. $V'$) be obtained from $W$ (resp. $V$) 
by removing the representatives of all cosets $Ax$, $x\in D$ (resp. $xD$, $x\in A$).
We can assume that $B'\subset W'$ and $C'\subset V'$, where $B=B'\cup\{e\}$ and $C=C'\cup\{e\}$.
We will have to perform at most 
\begin{equation}\label{subgroup_subgroup}
\binom{|G:A|-|D|}{|B|-1}\cdot\binom{|G:D|-|A|}{|C|-1}
\end{equation}
equality checks of $G=ABCD$.
For the sake of brevity, the process for computing the product of $ABCD$ will be referred to as the "$ABCD$ procedure" in the following.

This is the technique we used in Lemma \ref{lemma_A332252} where  
$A$ and $D$ are subgroups of $A_6$ and $|A|=9$, $|B|=|C|=2$, and $|D|=10$.
It follows from (\ref{subgroup_subgroup}) that 
we need to call the $ABCD$-procedure only $810$ times
and with fixed $A$ and $D$ the program written in \GAP 
found all $30$ factorizations in less than a second.

\textsc{Case 2.}
Suppose we have chosen a subgroup $A$ of $G$ and some subset $D\subset G$ such that 
the right cosets $Ax$, $x\in D$, are pairwise distinct, and $e\in D$ 
and $|D|$ is divisible by $|G|$.
We want to use a brute-force to compute sets $B$ and $C$ such that $G=ABCD$.
Let $W$ be a right transversal of $A$ in $G$ and 
$W'$ is obtained from $W$ by removing the representatives of right cosets $Ax$, $x\in D$. 
We will look for $B'$ such that $B'\subset W'$ and $B=B'\cup\{e\}$.
If $B$ is already selected, then we must choose $C'\subset V=G\setminus ABD$ and $C=C'\cup\{e\}$.
The number of possibilities to choose 
sets $B$ and $C$ in the specified way is 
\begin{equation}\label{subgroup_subset}
  \binom{|G:A|-|D|}{|B|-1}\cdot\binom{|G|-|A|\cdot|B|\cdot|D|}{|C|-1}.
\end{equation}
A similar deduction and formula are obtained if $D$ is a subgroup and $A$ is not a subgroup of $G$. 
In this case we should compute the factorization $G=DC^{-1}B^{-1}A^{-1}A^{-1}$ 
instead of the factorization $G=ABCD$.
The number of attempts is given by the formula
\begin{equation}\label{subset_subgroup}
  \binom{|G:D|-|A|}{|C|-1}\cdot\binom{|G|-|A|\cdot|C|\cdot|D|)}{|B|-1}.
\end{equation}

We used the method described here in lemmas \ref{lemma_6-7-2-2}, \ref{lemma_12-7-2}, and
\ref{lemma_A232352}.

In Lemma \ref{lemma_6-7-2-2} $G$ is the simple group of order $168$, $A$ is a subgroup of $G$ and $|A|=6$,
$D=\{e,d\}$ where $d$ has order $4$, and $|B|=7$, $|C|=2$.
Formula (\ref{subgroup_subset}) gives an upper bound on the number of calls to the $ABCD$-procedure:
$$
\binom{26}{6}\cdot\binom{84}{1}\approx19.3\cdot10^6.
$$
A program compiled on \GAP searched
through all the choices in about $10$ minutes and found $24$ solutions. The
solution given here is the fourth and the program found it in about $5$ seconds.

Two words about selecting a subset of $D$ after a subgroup of $A$ is fixed. 
Since $D=\{e,d\}$, by Corollary \ref{corollary_on two element factors} 
the order of $d$ is even and hence $2$ or $4$ in our case. 
However, if $d$ is of order $2$, then $D$ is a subgroup of order $2$. 
We know that all elements of order $2$ in the group $G$ are conjugate, 
so there are no factorizations of the form $AXD$ (see Lemma \ref{lemma_ATB}). 
On the other hand, we can check that 
any element of order $4$ of group $G$ can be taken as $d$. 

In the proof of Lemma \ref{lemma_12-7-2}, a subgroup $A$ of order $12$ is chosen 
in the group $G$ of order $168$. 
It is required to find only three factors such that $G=ABC$, 
where $|B|=7$, $|C|=2$. 
In accordance with the remark made above, we can take $C=\{e,c\}$, where $c$ has order $4$.
The number of calls to the $ABCD$-procedure in this case is equal to
$$
\binom{26}{6}=230230.
$$
A \texttt{GAP}-program found all $8$ solutions (with A and D fixed) in less than a second.

In Lemma \ref{lemma_A232352}, 
$D$ is a subgroup of order $10$, and the set $A=\{e,a\}\cdot\{e,t,t^2\}$ 
where $a$ and $t$ are elements of order $4$ and $3$, respectively (here the group $G=A_6$).
In addition $|B|=2$, $|C|=3$.
In view of (\ref{subset_subgroup}) the number of accesses to the $ABCD$-procedure is equal to 
$$
\binom{30}{2}\cdot\binom{180}{1}=78300.
$$
With $A$ and $D$ fixed, the \GAP program found $288$ solutions in $7$ seconds.
The first solution was given instantly.

In Lemma \ref{lemma_A18-522} we have:
$|G|=360$, $|A|=18$, $|D|=2$, $|B|=5$, $|C|=2$ by formula (\ref{subgroup_subset}) we obtain 
$$
\binom{18}{4}\cdot180=550800.
$$ 
With $A$ and $D$ fixed, the \GAP program found $144$ solutions in about $70$ seconds.
The first solution as in the case above was obtained instantaneously.

\textsc{Case 3.}
Suppose that a subgroup $A$ of the group $G$ and sets $B$, $D$ are chosen such
that 
\begin{equation}\label{choose_BD}
    |ABD|=|A|\cdot|B|\cdot|D|=m
\end{equation}
for some divisor $m$ of $|G|$.
Equality (\ref{choose_BD}) means, the right cosets of $Ax$ with representatives 
from $BD$ are pairwise distinct.
Let us show how we can find the set $C$.
First we construct a set $V$ as follows
\begin{equation}\label{vertices_gamma}
  V=\{x\in G\mid\ |ABxD|=m\}.
\end{equation}
Now consider a graph $\Gamma$ whose vertices are $V$ and
a pair of vertices $x\neq y$ connected by an edge if and only if
\begin{equation}\label{edges_gamma}
  ABxD\cap AByD\neq\varnothing.
\end{equation}
Suppose that we have found an independent set $C=\{c_1,\ldots,c_s\}$ of vertices of graph $\Gamma$.
Recall that an independent set of $\Gamma$ is a set of vertices of $\Gamma$,
no two of which are joined by an edge.
We want to check that each element $g\in ABCD$ has a single representation as 
$g=abcd$ with $a\in A$, $b\in B$, $c\in C$, $d\in D$.
In fact, if $abcd=a'b'c'd'$ and $a,a'\in A$, $b,b'\in B$, $c,c'\in C$ and $d,d'\in D$, then
$$
ABcD\cap ABc'D\neq\varnothing.
$$
If $c\neq c'$, then by virtue of (\ref{edges_gamma}) $c$ and $c'$ are joined by an edge in $\Gamma$.
But the latter is impossible since $c,c'\in C$ and $C$ is an independent set.
If $c=c'$ and for instance $a\neq a'$, then $|ABcD|<m$ but since $c\in C\subset V$
this contradicts condition (\ref{vertices_gamma}). 
The cases $b\neq b'$ and $d\neq d'$ are similarly treated.
If $s=\dfrac{|G|}{m}$, then we get $G=ABCD$, since in this case $|ABCD|=|G|$.

The outline of our algorithm for computing $C$ is as follows:
\begin{enumerate}
  \item Fix a subgroup $A$ of $G$, a right transversal $W$ of $A$ in $G$,
and a set $D$.
  \item For each  $B\subset W$ such that condition (\ref{choose_BD}) is satisfied do:
  \item By (\ref{vertices_gamma}) construct the set $V\subset G$. 
  \item Construct the graph $\Gamma$ with vertex set $V$ and edges specified in (\ref{edges_gamma}).
  \item Compute a maximal independent set $C$ of $\Gamma$.
  \item If $|C|=\dfrac{|G|}{m}$, then return $C$.
\end{enumerate}
This algorithm was implemented on \texttt{GAP} with the \texttt{GRAPE} package \cite{Soicher}.
We used the \texttt{IndependentSet} function of this package, which, as the author of the package write,
``returns a (hopefully large) independent set of the graph $\Gamma$''.
The main time cost in this computation comes from calling \texttt{IndependentSet} (point 5 of the algorithm). 
It is easy to see that the number of calls to point 5 is equal to
\begin{equation}\label{calling_IndependetSet}
  \binom{|G:A|-D}{|B|-1}.
\end{equation}

The described method of computation was applied in the search for factorizations 
in Lemmas \ref{lemma_6-2-7-2} and \ref{lemma_A223352-322352}.

In Lemma \ref{lemma_6-2-7-2} 
as a subgroup $A$ of the group $G$ of order $168$ we choose a subgroup of order $6$ and
$|D|=|B|=2$, $|C|=7$.
According to formula (\ref{calling_IndependetSet}) 
we will need to call the \texttt{IndependentSet} function $26$ times.
The first solution is found in half a second, and $21$ solutions are found in $10$ seconds.

In Lemma \ref{lemma_A223352-322352} 
as a subgroup $A$ of the group $A_6$ we choose a subgroup of order $12$ and
$|D|=2$, $|B|=3$, $|C|=5$.
We will need to call the \texttt{IndependentSet} function $378$ times.
The first solution is found in $10$ seconds and
$212$ solutions are found in about $10$ minutes.

\begin{remark}
In Lemma \ref{lemma_6-2-7-2}, we cannot take a subgroup of order $12$ as a subgroup of $A$ 
because every subgroup of order $12$ in $G$ is isomorphic to $A_4$ and $A_4$ has no $(2,3,2)$-factorizations. 

It would seem that we can search for sets $B$ and $C$ in the same way as in \textsc{Case 2}.
However, if we choose as $A$ a subgroup of order $6$ and as $D$ a set of two elements, 
we obtain $|B|=2$, $|C|=7$ and in view of (\ref{subgroup_subset}) the required number of calls of the $ABCD$-procedure is equal to
$$
\binom{26}{1}\cdot\binom{144}{6}>289\cdot10^9.
$$
As a result, the computational difficulties become insurmountable.

Similarly, if in constructing the factorization of $ABCD$ in Lemma \ref{lemma_A223352-322352} by the method outlined in \textsc{Case 2}, a subgroup of order $12$ is chosen as $A$ and $D$ consists of two elements, then the number of calls to the $ABCD$-procedure is equal to
$$
\binom{34}{2}\cdot\binom{324}{4}>252\cdot10^9
$$
and this is absolutely unacceptable.
\end{remark}

\section{Questions}\label{section-Questions}
It seems that the following must be true. 
\begin{conjecture}\label{conjecture_An}
The alternating group $A_n$ is multifold-factorizable whenever $n>5$ .
\end{conjecture} 
It may be easier to prove the following
\begin{conjecture} \label{conjecture_Sn}
The symmetric group $S_n$ is multifold-factorizable for all $n\geq2$.
\end{conjecture} 
By Lemma \ref{lemma_SimpleGeneral}\ref{normal_subgroup}, 
the positive solution of Conjecture \ref{conjecture_An} implies 
the positive solution of Conjecture \ref{conjecture_Sn}. 

The converse is not true.
For example, we know that $S_n$ is multifold-factorizable for $2\leq n\leq6$
 but $A_4$ and $A_5$ are not multifold-factorizable.

\end{document}